\documentclass[12pt]{article}

\usepackage{amsmath,amssymb,amsthm,amscd,dsfont}
\usepackage{amsfonts,latexsym,rawfonts}
\usepackage{lscape}
\usepackage{amscd, float,times,rotating}
\usepackage{pb-diagram}
\usepackage{hyperref}
\numberwithin{equation}{section}

\newtheorem{prop}{Proposition}[section]
\newtheorem{theorem}[prop]{Theorem}
\newtheorem{lemma}[prop]{Lemma}
\newtheorem{corollary}[prop]{Corollary}
\newtheorem{remark}[prop]{Remark}

\newtheorem{definition}[prop]{Definition}

\renewcommand{\d}{\mathrm{d}}

\newcommand{\loc}{\mathrm{loc}}

\newcommand{\Div}{\operatorname{div}}

\newcommand{\R}{\mathbb{R}}

\newcommand{\br}[1]{\left( #1 \right)}

\newcommand{\Ric}{\operatorname{Ric}}

\title{\bfseries Principal $p-$frequency estimates on non-compact manifolds with negative Ricci curvature}
\author{Xiaoshang Jin\thanks {X. J. is supported by ``the Fundamental Research Funds for the Central Universities'', HUST: \# 2025BRSXB002 and NSFC (Grant No. 12471054)}\ \ ,\ Zhiwei L\"u}
\date{}

\begin{document}

	\maketitle

    \begin{abstract}
   We establish a lower bound for the principal $p-$frequency $\lambda_{1,p}(\Omega)$ on a bounded domain $\Omega$ in a non-compact Riemannian manifold of dimension $n.$ Under the assumption that the Ricci curvature satisfies $\operatorname{Ric} \geq (n-1)K$ with $K<0,$ we prove that $\lambda_{1,p}(\Omega) > \bar{\lambda}_{D,K,n}$, where $D$ is the diameter of $\Omega$ and $\bar{\lambda}_{D,K,n}$ is explicitly defined as the first eigenvalue of an associated one-dimensional ordinary differential equation model that incorporates both $D$ and $K.$ Moreover, the estimate is sharp. This work extends previous results for the case $K=0$ to the geometrically more complex setting of negative Ricci curvature, and providing a new quantitative connection between the eigenvalue, the diameter of domains, and the curvature lower bound.
    \end{abstract}

    \section{Introduction}

    Let $p>1$ and $\Omega$ be a bounded domain with (piecewise) smooth boundary $\partial \Omega \neq \varnothing$ in a complete Riemannian manifold $(M,g).$  Define $\lambda_{1,p}(\Omega)$  to be the principal $p-$frequency of $\Omega$  if it is the first  positive eigenvalue of $\Omega$ with the Dirichlet condition:
        $$
        \left\{\begin{array}{ll}
            \Delta_p u = -\lambda_{1,p}(\Omega)|u|^{p-2}u & \quad \text{in $\Omega$}, \\
            u=0 & \quad \text{on $\partial \Omega$}.
        \end{array}\right.
    $$
    We call $u$ the eigenfunction of $\lambda_{1,p}(\Omega)$ with respect to the $p-$Laplace operator $\Delta_p:$
    $$
        \Delta_p v = \Div(|\nabla v|^{p-2}\nabla v), \quad v \in W^{1,p}_\loc(M)
    $$
    where the equality holds in the weak $W^{1,p}_\loc(M)$ sense. The principal $p-$frequency is also called the first Dirichlet eigenvalue of $\Delta_p,$ and has the following variational characterization:
    \begin{equation}\label{eq:variational}
        \lambda_{1,p}(\Omega) = \inf_{v \in W^{1,p}_0(\Omega) \setminus \{0\}} \frac{\int_\Omega |\nabla v|^p \,\d v_g}{\int_\Omega |v|^p \,\d v_g}.
    \end{equation}

    In this paper, we obtain a sharp estimate of $\lambda_{1,p}(\Omega)$  if the ambient manifold $(M,g)$ is non-compact and its Ricci curvature is bounded from below by a negative constant. To state our main theorem, we first introduce the corresponding one-dimensional ordinary differential equation model:
        \begin{definition}\label{def:barlambda}
        For arbitrary $D>0,\ K<0$ and $n\geq 2,$ let $\bar \lambda = \bar \lambda_{D,K,n}>0$ be the first eigenvalue of $[0,D]$ of the following eigenvalue problem:
        \begin{equation}\label{eq:1-dim-diri}
            \left\{\begin{array}{l}
            \frac{\d}{\d t} (\dot w^{(p-1)}) + (n-1)\sqrt{-K} \dot w^{(p-1)} + \bar\lambda w^{(p-1)} = 0, \\
            w(0) = 0, \dot w(D) = 0,
        \end{array}\right.
    \end{equation}
        where $$w^{(p-1)}:=|w|^{p-2}w.$$
        Proposition \ref{prop:2nd-model-existence} below shows that there is a smallest positive $\bar\lambda$ such that \eqref{eq:1-dim-diri} has a nontrivial solution which is unique in the sense of scaling.
    \end{definition}
    \begin{theorem}\label{thm:main-thm}
        Assume that $(M,g)$ is a complete and non-compact manifold of dimension $n$ whose Ricci curvature satisfies $\Ric\geq(n-1)K$ with $K<0.$ If $\Omega$ is a bounded smooth domain in $(M,g)$ with diameter $D,$ then we have the sharp estimate:
        \begin{equation}\label{eq:thm1.2}
          \lambda_{1,p}(\Omega) > \bar \lambda_{D,K,n}.
        \end{equation}
    \end{theorem}

    \begin{remark}\label{rmk:13}
    \begin{itemize}
      \item [(1)] By the variational method,
        \begin{equation}\label{eq:reyleigh-1-dim}
            \bar\lambda_{D,K,n} = \inf\left\{ \frac{\int_0^D |\dot u(t)|^p e^{(n-1)\sqrt{-K}t} \,\d t}{\int_0^D |u(t)|^p e^{(n-1)\sqrt{-K}t} \,\d t}:
            u\in W^{1,p}[0,D],\ u(0)=0,\ u\not\equiv 0 \right\}
        \end{equation}
        and the minimizer of this variational problem is the solution to \eqref{eq:1-dim-diri}.
      \item [(2)] If $K=0,$ then the equation \eqref{def:barlambda} is also solvable and
      $$\bar \lambda_{D,0,n}=(p-1)\Big(\frac{\pi_p}{2D}\Big)^p,\ \ w(t)=\sin_p\br{\frac{\pi_p}{2D}t}$$
          where $\sin_p$ and $\pi_p$ are the generalized trigonometric functions (see Definition \ref{def:gen-tri}). In particular, for $p=2,$
           $$\bar \lambda_{D,0,n}=\Big(\frac{\pi}{2D}\Big)^2,\ \ w(t)=\sin\br{\frac{\pi}{2D}t}.$$
           For this special case, Theorem 1.2 was already proved in \cite{wang2022sharp} for $p=2$ and in \cite{jin2025estimate} for general $p>1.$
      \item [(3)] The lower bound $\bar\lambda_{D,K,n}$ in \eqref{eq:thm1.2} is never attained on any smooth bounded domain. Nevertheless, this estimate is sharp in the following sense: one can construct a sequence of bounded domain $\Omega_i\subseteq (M_i,g_i)$ with diameter $D_i$ such that
          $$\lim\limits_{i\rightarrow\infty}D_i=D\ \ \text{and}\ \ \lim\limits_{i\rightarrow\infty}\lambda_{1,p}(\Omega_i) = \bar \lambda_{D,K,n}.$$
    \end{itemize}
    \end{remark}
    To better understand the geometric behavior of the model eigenvalue $\bar{\lambda}_{D,K,n}$, we examine its asymptotic expansions in the two extremal regimes: very small diameters and very large diameters.
    \begin{theorem}\label{thm2}
    Let $\bar\lambda_{D,K,n}$ be defined as in Definition \ref{def:barlambda}, then
      \begin{equation}\label{eq:15}
        \bar\lambda_{D,K,n}=(p-1)\br{\frac{\pi_p}{2D}}^p+O(\frac{1}{D^{p-1}}),\ \ \text{as}\ D\rightarrow 0,
        \end{equation}
        \begin{equation}\label{eq:16}
       \ln \bar\lambda_{D,K,n}=-(n-1)\sqrt{-K}\cdot D+O(1),\ \ \text{as}\ D\rightarrow\infty.
      \end{equation}
      where $O(1)$ denotes a quantity bounded independent of $D.$
    \end{theorem}
    The first formula shows that for tiny domains the curvature lower bound $K$ becomes negligible; the eigenvalue behaves exactly as in the zero curvature case, i.e. like the first Dirichlet eigenvalue of the $p-$Laplacian. The correction $O(\frac{1}{D^{p-1}})$ captures the slight influence of the curvature term in the model equation. As $D\to\infty$, the eigenvalue decays exponentially. This reflects how negative curvature lets the eigenfunction extend with minimal slope, sharply lowering the energy. Furthermore, if $p=2,$ then the exact formula of $\lambda_{D,K,n}$ could be achieved according to \eqref{eq:relation}.
    \par Combining the two Theorems above, we derive that
    \begin{corollary}\label{coro1.5}
      Let $(M,g)$ be a complete noncompact $n$-manifold with $\mathrm{Ric}\ge (n-1)K$, $K<0$. Let $\Omega\subset M$ be a bounded smooth domain of diameter $D=\operatorname{diam}(\Omega)$.
      \begin{itemize}
\item For any fixed $D_0>0$, there exists a constant
$C_1=C_1(p,n,K,D_0)>0$ such that
$$
\lambda_{1,p}(\Omega)\ge \frac{C_1}{D^p}\qquad \text{for all }0<D\le D_0.
$$

\item There exists a constant $C_2=C_2(p,n,K)>0$ independent of $D$ such that
$$\lambda_{1,p}(\Omega)\ge C_2\,e^{-(n-1)\sqrt{-K}\,D}\qquad \text{for all }D>0.
$$
\end{itemize}
    \end{corollary}
Moreover, one can show that $C_1$ is decreasing of $D_0$ and tends to $0$ as $D_0$ approaches infinity.
\\~ 
    \par A classical approach for estimating the first Dirichlet eigenvalue (in the linear case $p=2$) is the gradient estimate for eigenfunctions introduced by Li--Yau \cite{li1980estimates}. This method can be generalized to the case $p>1$. For further studies on eigenvalue problems of the $p-$Laplace operator we refer to \cite{JIN2025110603} \cite{kawai2003first}\cite{matei2000first}\cite{miao2025first}\cite{naber2014sharp}\cite{valtorta2012sharp}\cite{zhang2007lower}.

    Recall that in \cite{naber2014sharp}, the authors found a sharp estimate for the Neumann $p$-eigenvalue via gradient estimates and a careful analysis of one dimensional models, while \cite{wang2022sharp}, the maximum principle was applied to a function built from the Busemann function.   Inspired by these works, we employ a Barta's inequality with a test function defined by a one-dimensional model to derive the sharp lower bound.
    \\~\par Here is the outline of this paper. In Section 2 we provide some preliminary materials concerning the generalized trigonometric functions.  Subsequently, we will prove  the existence of a smallest positive $\bar\lambda$ such that the equation \eqref{eq:1-dim-diri} admits a unique solution satisfying  $w(D)=1$ and prove Theorem \ref{thm2} and Corollary \ref{coro1.5}.
    \par In Section 3, we present the proof of Theorem 1.2. First, by utilizing the Busemann function and the solution to  \eqref{eq:1-dim-diri}, we construct an appropriate test function. The desired lower bound for the principal $p-$frequency is then derived by applying a Barta-type inequality. Second, to demonstrate the sharpness of this estimate, we construct a sequence of domains within warped product manifolds satisfying $\Ric \geq (n-1)K$, and show that the limit of their principal $p-$frequencies coincides exactly with $\bar{\lambda}_{K,D,n}.$

    \section{Analysis of the one dimensional model}

    First we introduce the generalized trigonometric functions. For the following conclusions about them, one may refer to \cite{dosly2005half}.
    \begin{definition}\label{def:gen-tri}
        Define
    $$
        \arcsin_p(x) := \int_0^x \frac{1}{(1-t^p)^{1/p}} \,\d t, \quad x \in [0,1]
    $$
    and
    $$
        \pi_p := 2\arcsin_p(1) = \frac{2 \pi}{p \sin(\pi/p)}.
    $$
    Let $\sin_p(t)$ be the inverse function of $\arcsin_p$ for $t \in [0,\pi_p/2]$ and set
    $$
        \sin_p(t) = \sin_p(\pi_p-t), \quad t \in [\pi_p/2,\pi_p].
    $$
    Then we extend $\sin_p$ to be an odd $2\pi_p$-periodic function defined on the whole line $\R$.  Let $\cos_p(t) := (\sin_p)'(t)$, then we have
    $$
        |\sin_p(t)|^p + |\cos_p(t)|^p = 1, \quad \forall t \in \R.
    $$
    It is easy to see that $\sin_p,(\cos_p)^{p-1} \in C^1(\R)$.
    \end{definition}

    Now introduce the one-dimensional model.

    \begin{definition}
        For any $\bar\lambda>0$. Let $w$ be the solution of the following boundary-value problem on $\R$:
    \begin{equation}\label{eq:1-dim-model}
        \left\{\begin{array}{l}
            \frac{\d}{\d t} (\dot w^{(p-1)}) + (n-1)\sqrt{-K} \dot w^{(p-1)} + \bar\lambda w^{(p-1)} = 0, \\
            w(0) = 0, w(D) = 1.
        \end{array}\right.
    \end{equation}
    \end{definition}

    \begin{remark}\label{rmk:1}
        Let $M = [0,D] \times_{e^t} \mathbb{S}^{n-1}$. Then $u(x)=w(r(o,x))$ solves $\Delta_p u=- \bar\lambda u^{(p-1)}.$
    \end{remark}

    We now show that for arbitrary $D>0,\ K<0,$ there is a unique $\bar \lambda=\bar \lambda_{D,K,n}>0$ such that \eqref{eq:1-dim-model} has a solution with the additional condition $\dot w(D)=0.$ Apply the Pr\"ufer transformation to \eqref{eq:1-dim-model}:
    $$
        \alpha w = e \sin_p \phi, \quad \dot w = e \cos_p \phi,
    $$
    where $e>0$ and $\alpha=\br{\frac{\bar\lambda}{p-1}}^{1/p}$. Then $\phi$ satisfies the following equation
    \begin{equation}\label{eq:prufer-phi}
        \left\{\begin{array}{l}
            \dot\phi = \alpha + \frac{(n-1)\sqrt{-K}}{p-1} \cos_p^{p-1} \phi \sin_p \phi =: \alpha+F(\phi), \\
            \phi(0) = 0 \mod \pi_p.
        \end{array}\right.
    \end{equation}
    Moreover, $\dot w(D)=0$ is equivalent to $\phi(D) = \pi_p/2 \mod \pi_p$.

    \begin{lemma}\label{prop:2nd-model-existence}
        There exists a unique $\alpha>0$ such that \eqref{eq:prufer-phi} has a solution $\phi$ satisfying that $\phi(D)=\pi_p/2.$ Consequently, there exists a smallest $\bar\lambda>0$ such that \eqref{eq:1-dim-model} has a non-trivial solution $w$ satisfying that $\dot w(D)=0$.
    \end{lemma}

    \begin{proof}
        Let $\phi$ be a solution of \eqref{eq:prufer-phi} such that $\phi(0)=0$, such $\phi$ exists by the unique existence theorem of ODEs. Denote by $\delta(\alpha) > 0$ the smallest number such that $\phi(\delta) = \pi_p/2$ (let $\delta=+\infty$ if such a number does not exist).

        If $\phi \in [0,\pi_p/2]$, we have $F(\phi) = \frac{(n-1)\sqrt{-K}}{p-1} \cos_p^{p-1} \phi \sin_p \phi\geq 0.$ Hence $\dot \phi = \alpha + F(\phi) > 0$, and $\phi$ is strictly increasing. Thus we have
        $$
            \delta = \int_0^{\delta} \frac{\dot\phi}{\alpha+F(\phi)} \d t = \int_0^{\pi_p/2} \frac{\d\phi}{\alpha+F(\phi)}.
        $$
        Therefore $\delta \to 0$ as $\alpha \to \infty$. Note that $F(\phi) \sim \sin_p \phi$ as $\phi \to 0$, we have $|F(\phi)| \leq C\phi$ in $[0,\varepsilon]$ for some small $\varepsilon>0.$ Hence
        $$
            \delta \geq \int_0^{\pi_p/2} \frac{\d\phi}{\alpha+F(\phi)} \geq \int_0^\varepsilon \frac{\d\phi}{\alpha+C\phi}.
        $$
        Therefore $\delta \to \infty$ as $\alpha \to 0$. Since $\delta$  depends continuously on $\alpha$, for arbitrary $D>0$, there is a unique $\alpha$ such that $\delta=D$.
    \end{proof}

Lemma \ref{prop:2nd-model-existence} indicates that for arbitrary $D>0$, there is a smallest $\bar \lambda > 0$ such that \eqref{eq:1-dim-diri} has a solution in the sense of scaling, and the solution is non-negative on $[0,D]$.
    \begin{remark}\label{rmk25}
    \begin{itemize}
      \item [(1)] From the above proof we obtain the relation
        \begin{equation}\label{eq:relation}
          D = \int_0^{\pi_p/2} \frac{\d \phi}{\br{\frac{\bar\lambda}{p-1}}^{1/p} + \frac{(n-1)\sqrt{-K}}{p-1} \cos_p^{p-1}\phi \sin_p \phi}.
        \end{equation}
        Therefore for any $c>0,$
        $$
           \bar\lambda_{c^{-1}D,c^2 K, n}=c^{p} \bar\lambda_{D,K,n}.
        $$
      \item [(2)] Let $\Omega$ be a bounded domain in $(M,g)$ and $c$ is a positive constant. If
   $$\begin{cases}
   \text{diam}(\Omega,g)=D\\
   \Ric[g]\geq (n-1)K,
    \end{cases}
   $$
   then the scaling metric $(M,c^{-2}g)$ satisfies that
   $$\begin{cases}
        \text{diam}(\Omega,c^{-2}g)=c^{-1}D \\
       \Ric[c^{-2}g]\geq (n-1)c^2K\\
      \lambda_{1,p}(\Omega,c^{-2}g)=c^{p}\lambda_{1,p}(\Omega,g).
     \end{cases}
   $$
   \item [(3)] As a consequence of (1) and (2), we can derive that the estimate \ref{eq:thm1.2} is scaling-invariant.
    \end{itemize}
    \end{remark}
    \begin{proof}[Proof of Theorem \ref{thm2}]
    From the proof of Lemma \ref{prop:2nd-model-existence} or \eqref{eq:relation}, we know that $D$ and $\bar\lambda$ are in one-to-one correspondence. Furthermore,
    $$
    D\rightarrow 0\Leftrightarrow\bar\lambda\rightarrow\infty,\ \ \ \ D\rightarrow \infty\Leftrightarrow\bar\lambda\rightarrow0.
    $$
    We rewrite \eqref{eq:relation} as
    $$
    \begin{cases}
       D = \int_0^{\frac{\pi_p}{2}} \frac{\d \phi}{\alpha + \beta \cos_p^{p-1}\phi \sin_p \phi}\\
      \alpha=\br{\frac{\bar\lambda}{p-1}}^{1/p} \\
      \beta=\frac{(n-1)\sqrt{-K}}{p-1}.
    \end{cases}
    $$
    \begin{itemize}
      \item [(1)] If $\alpha\rightarrow\infty,$ then
      $$\alpha D = \int_0^{\frac{\pi_p}{2}} \frac{\d \phi}{1 + \frac{\beta}{\alpha} \cos_p^{p-1}\phi \sin_p \phi}=\frac{\pi_p}{2}+O(\frac{1}{\alpha})$$
      As a consequence,
      $$\alpha=\frac{\pi_p}{2D}+O(1), \ \ D\rightarrow 0
      $$
      which would imply \eqref{eq:15}.
      \item [(2)] If $\alpha\rightarrow 0,$ then we first notice that
      $$
      \cos_p^{p-1}\phi \sin_p\phi=
      \begin{cases}
          \phi+O(\phi^{p+1})\ &\ \text{as}\ \phi\rightarrow 0,\\
          (p-1)(\frac{\pi_p}2-\phi)+O((\frac{\pi_p}2-\phi)^{1+\frac{p}{p-1}})&\ \text{as}\ \phi\rightarrow \frac{\pi_p}2.
        \end{cases}
      $$
      Hence we could define the funciton
      $$
      h(\phi)=
      \begin{cases}
         \frac{1}{\alpha + \beta \phi}\ &\ \text{if}\ \phi\in[0,\frac{\pi_p}4)\\
          \frac{1}{\alpha + \beta(p-1)(\frac{\pi_p}2-\phi)} & \text{if}\ \phi\in[\frac{\pi_p}4,\frac{\pi_p}2].
        \end{cases}
      $$
      It is clear that
      $$
      \int_0^{\frac{\pi_p}{2}} \frac{\d \phi}{\alpha + \beta \cos_p^{p-1}\phi \sin_p \phi}-\int_0^{\frac{\pi_p}{2}} h(\phi)\d \phi
      $$
      is bounded and
      $$\begin{aligned}
      \int_0^{\frac{\pi_p}{2}} h(\phi)\d \phi & =\int_0^{\frac{\pi_p}{4}}\frac{1}{\alpha + \beta \phi}\d \phi+\int_{\frac{\pi_p}{4}}^{\frac{\pi_p}{2}} \frac{1}{\alpha + \beta(p-1)(\frac{\pi_p}2-\phi)}\d \phi
      \\&=\frac1\beta\ln\frac{\alpha+\beta\frac{\pi_p}{4}}{\alpha}+\frac1{\beta(p-1)}\ln\frac{\alpha+\beta(p-1)\frac{\pi_p}{4}}{\alpha}
      \\&=\Big(\frac1\beta+\frac1{\beta(p-1)}\Big)\ln\frac1\alpha+O(1)
      \end{aligned}
      $$
      Finally, we derive that
      $$
      D=\Big(\frac1\beta+\frac1{\beta(p-1)}\Big)\ln\frac1\alpha+O(1)=\frac{-\ln\bar\lambda}{(n-1)\sqrt{-K}}+O(1)
      $$which would imply \eqref{eq:16}.
    \end{itemize}
    \end{proof}
    
    \begin{proof}[Proof of Corollary \ref{coro1.5}]
    \begin{itemize}
      \item For any fixed $D_0>0$ and $\Omega\subseteq M$ with ${\rm diam}(\Omega)=D<D_0,$ we only need to show that
      $$\bar\lambda_{D,K,n}\geq \frac{C_1}{D^p}.$$
      \eqref{eq:15} implies that 
      $$\bar\lambda_{D,K,n} D^p=(p-1)\br{\frac{\pi_p}{2}}^p+O(D),\ \ \text{as}\ D\rightarrow 0.$$
      According to \eqref{eq:relation}, $\bar\lambda$ is a continuous function of $D$ on $(0,D_0],$ hence
      $$\bar\lambda_{D,K,n} D^p\geq C_1(p,n,K,D_0)>0\ \ {\rm on}\ \ (0,D_0].$$
      \item For the second case, by \eqref{eq:16},
      $$
      \ln \bar\lambda_{D,K,n}+(n-1)\sqrt{-K}\cdot D\  \text{is bounded on}\ [1,+\infty)
      $$ and it implies that
      \begin{equation}\label{eq:2.4}
      \bar\lambda_{D,K,n}e^{(n-1)\sqrt{-K}D}\ge C_2'(p,n,K)>0\  \text{for}\ \ D\in[1,+\infty). 
      \end{equation}
        Then we finish the proof as \eqref{eq:2.4} is obvious true for $D\in(0,1].$
    \end{itemize}
    \end{proof}
    \par We end this section by quoting a result from \cite{jin2024lower} that will be essential later.

    \begin{theorem}\label{thm:1.1}
        Let $(M,g)$ be a complete Riemannian manifold and $\Omega \subseteq M$ a bounded domain. If there exists $f \in W^{1,p}_\loc(\Omega)$ and a constant $\mu > 0$ such that
        $$
            \Delta_p f - (p-1) |\nabla f|^p \geq \mu
        $$
        in the weak $W^{1,p}_\loc$ sense, then
        $$
            \lambda_{1,p}(\Omega) \geq \mu.
        $$
        Moreover, if $\lambda_{1,p}(\Omega) = \mu$, then $e^{-f}$ is a $p$-eigenfunction corresponding to $\lambda_{1,p}(\Omega)$.
    \end{theorem}

    A complete proof can be found in \cite{jin2025estimate}.

    \section{Proof of the main theorem}
    The proof of Theorem \ref{thm:main-thm} is divided into two parts: establishing the lower bound and proving its sharpness.
   \subsection{The lower bound}
    \par Since $M$ is non-compact, we can select a ray $\gamma$ on $M$ and define the Busemann function
    $$\beta(x) = \beta_\gamma(x) = \lim_{t \to \infty} (d(x,\gamma(t))-t).$$ Then by the Laplacian comparison,
    $$
        \Delta \beta \leq (n-1)\sqrt{-K}\ \text{on}\ M
    $$
    as $\Ric\geq(n-1)K.$ Moreover, $\beta$ is a Lipschitz continuous function and
    $$
    |\nabla \beta|=1 \ \text{almost everywhere in $M$}
    $$
   (see \cite{busemann2005geometry}). Therefore, if the diameter of $\Omega$ is $D,$ then there is a constant $A\in\R$ such that
   $$\beta(x)\in [A,A+D],\ \ \forall x\in\Omega.$$
   Without loss of generality we may assume $A=0,$ otherwise we replace $\beta$ by $\beta-A.$
    \par Recall that $w$ is the unique positive solution of equation \eqref{eq:1-dim-diri} with $\bar\lambda=\bar\lambda_{K,D,n}$ given by Lemma \ref{prop:2nd-model-existence}, normalized so that s $w(D)=1.$ Let us  consider the function
    $$
        f(t) = -\ln w(t),\ \ t\in(0,D].
    $$
    Therefore,
    $$
        \dot f(t) = -\frac{\dot w(t)}{w(t)}.
    $$
    We claim that $\dot f(t) \leq 0,$ or equivalently, $w$ is an increasing function in $[0,D].$  Indeed, at any critical points of $w$, we have
     $$(p-1)|\dot w|^{p-2}\ddot w = (\dot w^{(p-1)})'=-\bar \lambda w < 0$$ when $t \in (0,D]$.
     Therefore every critical point of $t$ in $(0,D]$ must be a local maximal point, and $t=D$ is such a point.
      Suppose $\dot w(t_0)<0$ for some
     $t_0\in (0,D),$ then $\dot w<0$ in a neighbourhood of $t_0$ since $w\in C^{1,\alpha}[0,D]$ according to the regularity theory. Then we can assume that $w$ decreases in $[a,b]$ ($b \neq D$ since $D$ is a maximal point), and does not decrease in $[a,b+\varepsilon]$, $\forall \varepsilon > 0$. Then $b$ is a local minimal point, which contradicts to our conclusion that $w$ has only local maximal points. Therefore $\dot w(t) \geq 0$. Moreover, we can get that $w$ is strictly increasing function.

    Let us compute
    \begin{align*}
        (\dot f^{(p-1)})' &= -\frac{(\dot w^{(p-1)})'w^{(p-1)}-\dot w^{(p-1)}(w^{(p-1)})'}{w^{2(p-1)}} \\
        &= -\frac{-(n-1)\sqrt{-K} \dot w^{(p-1)}-\bar\lambda w^{(p-1)}}{w^{(p-1)}} + (p-1)\frac{|\dot w|^p |w|^{p-2}}{w^{2(p-1)}} \\
        &= -(n-1)\sqrt{-K} \dot f^{(p-1)} + \bar\lambda + (p-1) |\dot f|^p.
    \end{align*}
    Consider the test function
    $$
    f(\beta(\cdot)) \ \text{in}\ \Omega
    $$
    and a direct calculation indicates that
    $$
    \nabla f(\beta)=\dot{f}\nabla\beta
    $$ and
    \begin{equation}\label{eq:estimating}
    \begin{aligned}
        \Delta_p (f(\beta)) &= \Div (|\nabla (f(\beta))|^{p-2}\nabla(f(\beta))) \\
        &= \Div (\dot f^{(p-1)} \nabla \beta) \\
        &= (\dot f^{(p-1)})' + \dot f^{(p-1)} \Delta \beta \\
        &= (-(n-1)\sqrt{-K} + \Delta \beta) \dot f^{(p-1)} + \bar\lambda + (p-1)|\dot f|^p \\
        &\geq \bar\lambda + (p-1)|\dot f|^p.
    \end{aligned}
    \end{equation}
    Therefore
    $$
        \Delta_p(f(\beta)) - (p-1)|\nabla f(\beta)|^p \geq \bar\lambda.
    $$
    By Theorem \ref{thm:1.1},
    $$
        \lambda_{1,p}(\Omega) \geq \bar \lambda.
    $$

    Finally, we will show that the lower bound cannot be achieved. In fact, if $\lambda_{1,p}(\Omega) = \bar\lambda$ for some smooth bounded domain $\Omega,$ then the equality in the last line of \eqref{eq:estimating} is achieved, i.e.
    $$\Delta \beta = (n-1)\sqrt{-K}$$
     in the sense of distrubution.
     Moreover, by Theorem \ref{thm:1.1}, $e^{-f} = w(\beta(x))$ is a Dirichlet $p$-eigenfunction. Since $w$ is strictly increasing, we have that $$w(\beta(x))|_{\partial \Omega} \equiv 0\Rightarrow \beta(x)|_{\partial \Omega}=0.$$
     By the elliptic PDE theory, $\beta$ is a smooth function. Then
    $$|\nabla \beta|=1 \ \text{almost everywhere}$$
    would imply that $|\nabla \beta|=1$ everywhere. As a consequence, $\beta$ cannot take its maximum or minimum value in the interior of $\Omega.$ Then $\beta(x)\equiv 0$ in $\Omega$ as $\beta(x)|_{\partial \Omega}=0.$ Hence $\nabla \beta=0$ in $\Omega$ which is obviously a contradiction. Therefore $\lambda_{1,p}(\Omega) > \bar\lambda$.

   \subsection{Sharpness of the estimate}
   \par Thanks to Remark \ref{rmk25}, we only need to consider the case $K=-1.$
    \par To establish the sharpness of our lower bound, we  construct a family of warped product manifolds
$$
(M_\varepsilon,g_\varepsilon)=\bigl([0,\infty)\times\mathbb{S}^{n-1},\ dt^2+f_\varepsilon(t)^2\, g_{\mathbb{S}^{n-1}}\bigr).
$$
The warping function is designed so that $f_\epsilon(t)=\sinh t$ near the origin, ensuring smoothness, and transitions to the form $f_\epsilon(t)=\delta(\epsilon)\, e^{-t}$ for $t\geq\epsilon$. After the change of variable $s=t+\ln\delta(\varepsilon),$ the metric asymptotically resembles the cusp-type metric
$
ds^2+e^{-2s}g_{\mathbb{S}^{n-1}}.
$
In essence, we are taking a cylindrical segment of approximate length $D$ from the far end of this cusp (where
$s$ is large) and smoothly capping it off at one end by connecting to a hyperbolic ball.
\par Lemma \ref{lem:ricci-lower} shows that these manifolds satisfy
$\Ric\geq-(n-1).$
The geometry is highly anisotropic: the spherical directions shrink exponentially, so that geodesic balls of diameter about $D$ behave essentially like one-dimensional intervals. This model effectively captures the extremal behavior of the principal $p$-frequency under the given curvature constraint, and will be used to construct a sequence of domains whose eigenvalues converge to the theoretical lower bound $\bar{\lambda}_{D,K,n}.$
\\~\par  More concretely, fix $\varepsilon>0$ sufficiently small and let $\psi_\varepsilon$ be a smooth cutoff function such that
    $$
        \psi_\varepsilon(x) = \left\{\begin{array}{ll}
            0, & \quad 0 \leq x \leq \varepsilon/2, \\
            1, & \quad x \geq \varepsilon,
        \end{array}\right.
    $$
    and let
    $$
        h_\varepsilon(x) = (1-\psi_\varepsilon(x))\coth x - \psi_\varepsilon(x).
    $$
    Let
    $$
        f_\varepsilon(x) = \left\{\begin{array}{ll}
            \sinh x, & \quad 0 \leq x \leq \varepsilon/2, \\
            \sinh (\varepsilon/2) \exp \int_{\varepsilon/2}^x h_\varepsilon(t) \,\d t, & \quad x \geq \varepsilon/2.
        \end{array}\right.
    $$
    Then $f \in C^\infty[0,+\infty)$ and
    \begin{equation}\label{eq:f}
      f_\varepsilon(x) = \left\{\begin{array}{ll}
            \sinh x, & \quad 0 \leq x \leq \varepsilon/2, \\
            \delta(\varepsilon) e^{-x}, & \quad x \geq \varepsilon,
        \end{array}\right.
    \end{equation}
    where
    \begin{equation}\label{eq:delta}
    \begin{aligned}
        \delta(\varepsilon) & = \sinh(\varepsilon/2) \exp\left(\int^\varepsilon_{\varepsilon/2} h_\varepsilon(t)\,\d t + \varepsilon\right)
         \\ & \leq  \sinh(\varepsilon/2)\exp\left(\int^\varepsilon_{\varepsilon/2}\coth t \,\d t + \varepsilon\right)
         \\& = e^\varepsilon\sinh\varepsilon.
        \end{aligned}
    \end{equation}
    Now Consider a sequence of manifolds
    $$
        (M_\varepsilon,g_\varepsilon)=([0,+\infty) \times \mathbb{S}^{n-1},dt^2+f_\varepsilon(t)^2g_{\mathbb{S}^{n-1}}).
    $$
    Then the warped product manifolds are smooth since $f_\varepsilon(t)=\sinh t$ near $0.$ Moreover,
     \begin{lemma}\label{lem:ricci-lower}
        The manifolds $(M_\varepsilon,g_\varepsilon)$ satisfy that $\Ric[g_\varepsilon] \geq -(n-1)$.
    \end{lemma}

    \begin{proof}
       We omit the $\varepsilon$'s in subscripts for simplicity of notations and need to prove
    $$
        (n-2)\frac{1-\dot f^2}{f^2} - \frac{\ddot f}{f} \geq -(n-1).
    $$
    The inequality above holds on $x \in [0,\varepsilon/2] \cup [\varepsilon,+\infty]$. So we focus on the situation that $x \in [\varepsilon/2,\varepsilon]$. Note that $h = f'/f$, then the inequality is equivalent to
    \begin{equation}\label{ineq:ricci}
         h' + (n-1)( h^2 - 1) \leq \frac{n-2}{f^2}.
    \end{equation}
    Since $h(t) \leq \coth t$, we have
    $$\begin{aligned}
        f(x) & =\sinh (\varepsilon/2) \exp \int_{\varepsilon/2}^x h_\varepsilon(t) \,\d t
        \\ & \leq \sinh (\varepsilon/2)\exp\int_{\varepsilon/2}^x \coth t \,\d t
        \\ &=\sinh x.
        \end{aligned}
    $$
    Therefore, to prove \eqref{ineq:ricci}, it suffices to prove
    \begin{equation}\label{ineq:ricci1}
        h'(t) + (n-1)(h^2(t) - 1) \leq \frac{n-2}{\sinh^2 t} = \frac{n-1}{\sinh^2 t} + (\coth t)',
    \end{equation}
    which is equivalent to
    $$
         h'(t) - (\coth t)' + (n-1)(h^2 - \coth^2 t) \leq 0.
    $$
    Denote $g(t) = \coth t - h(t) \geq 0$, then \eqref{ineq:ricci1} is equivalent to
    \begin{equation}\label{ineq:ricci2}
         g'(t) + (n-1)g(t)(\coth t + h(t)) \geq 0.
    \end{equation}
    By the definition of $g$,
    \begin{align*}
        g(t) &= \psi(t)(\coth t + 1), \\
        g'(t) &= \psi'(t) (\coth t+1) - \frac{\psi(t)}{\sinh^2 t}.
    \end{align*}
    Insert into the LHS of \eqref{ineq:ricci2}, and using that $h(t)\geq -1,$ we get
    \begin{align*}
        \text{LHS of \eqref{ineq:ricci2}} &= \psi'(t)(\coth t+ 1)- \frac{\psi(t)}{\sinh^2 t} + (n-1)\psi(t)(\coth t + 1)(\coth t + h(t))
        \\ &\geq 0-\frac{\psi(t)}{\sinh^2 t} + (n-1)\psi(t)(\coth t + 1)(\coth t -1) \\
        & =\frac{(n-2) \psi(t)}{\sinh^2 t}\\
        &\geq 0.
    \end{align*}
    Therefore \eqref{ineq:ricci2} is proved, which implies that $\Ric[g_\varepsilon] \geq -(n-1).$
    \end{proof}
    With the preparations above, we can now prove the sharpness of the estimate \eqref{eq:thm1.2}:
    \begin{proof}[Proof of the sharpness]
    Fix a sufficiently small $\varepsilon>0$ and let $o$ be the center of $(M_\varepsilon,g_\varepsilon).$ Assume that
    $$\delta=\delta(\varepsilon)< e^\varepsilon\sinh\varepsilon$$ as defined in \eqref{eq:f} and select a point $z= \exp_o((1-\delta)\theta)$ for some $\theta \in \mathbb{S}^{n-1}.$ We claim that
    \begin{equation}\label{eq:BsubsetB}
      B(o,2-\delta-\pi \delta e^{\delta-1}) \subseteq B(z,1).
    \end{equation}
    Indeed, for any $x \in B(o,2-\delta-\pi \delta e^{\delta-1}),$ we have that
     $$
     x = \exp_o(r_x \theta_x)
     $$ for some $\theta_x \in \mathbb{S}^{n-1}$ and $r_x \leq 2-\delta-\pi \delta e^{\delta-1}.$
     Let $y_x = \exp_o((1-\delta)\theta_x)$, then
    $$\begin{aligned}
        d(x,z) &\leq d(x,\exp_o((1-\delta)\theta_x)) + d(\exp_o((1-\delta)\theta_x),z)
         \\ & \leq |r_x-(1-\delta)|+\pi f_\varepsilon(1-\delta)
         \\ & \leq 2-\delta-\pi \delta e^{\delta-1}-(1-\delta) + \pi \delta e^{\delta-1}
         \\ & = 1.
         \end{aligned}
    $$
    We consider the bounded domain
    $$
    \Omega_\varepsilon=B(o,2-\delta-\pi \delta e^{\delta-1})\subseteq M_\varepsilon
    $$
    with diameter $D_\varepsilon=\text{diam}(\Omega_\varepsilon).$ Then by equation \eqref{eq:BsubsetB},
    \begin{equation}\label{eq:Diameter}
      2-\delta-\pi \delta e^{\delta-1}\leq D_\varepsilon\leq 2
    \end{equation}
    In order to make an estimate of $\lambda_{1,p}(\Omega_\varepsilon),$ we consider the test function in $\Omega_\varepsilon:$
    $$\begin{cases}
       v(x)=w(R_\varepsilon-r(x)) \\
       R_\varepsilon=2-\delta-\pi \delta e^{\delta-1} \\
       r(x)=\text{distance function of $o$}\\
       w=\text{positive solution to \eqref{eq:1-dim-diri} when $D=R_\varepsilon$.}
      \end{cases}$$
    Then $v\in W^{1,p}_0(\Omega_\varepsilon)$ and hence
    \begin{align*}
        \lambda_{1,p}(\Omega_\varepsilon) &\leq \frac{\int_{\Omega_\varepsilon} |\nabla v(x)|^p \,\d v}{\int_{\Omega_\varepsilon} |v(x)|^p \,\d v} \\
        &= \frac{\int_0^{R_\varepsilon} |\dot w(R_\varepsilon-r)|^p \cdot \omega_{n-1}f^{n-1}_\varepsilon(r) \,\d r}{\int_0^{R_\varepsilon} |w(R_\varepsilon-r)|^p\cdot \omega_{n-1}f^{n-1}_\varepsilon(r) \,\d r} \\
        &= \frac{C_1(\varepsilon)+\int_\varepsilon^{R_\varepsilon} |\dot w(R_\varepsilon-r)|^p \cdot (\delta e^{-r})^{n-1} \,\d r}{C_2(\varepsilon)+\int_\varepsilon^{R_\varepsilon} |w(R_\varepsilon-r)|^p\cdot (\delta e^{-r})^{n-1} \,\d r}
        \\
        &= \frac{C_1(\varepsilon)-\int_0^{R_\varepsilon-\varepsilon} |\dot w(t)|^p \cdot (\delta e^{t-R_\varepsilon})^{n-1} \,\d t}{C_2(\varepsilon)-\int_0^{R_\varepsilon-\varepsilon} |w(t)|^p\cdot (\delta e^{t-R_\varepsilon})^{n-1} \,\d t}
    \end{align*}
    where
     $$\begin{cases}
        C_1(\varepsilon) = \int_0^{\varepsilon} |\dot w(R_\varepsilon-r)|^p \cdot \omega_{n-1}f^{n-1}_\varepsilon(r) \,\d r \\
        C_2(\varepsilon) = \int_0^{\varepsilon} | w(R_\varepsilon-r)|^p \cdot \omega_{n-1}f^{n-1}_\varepsilon(r) \,\d r
       \end{cases}\to 0$$
      as $\varepsilon \to 0$. Consequently,
    \begin{equation}\label{eq:9}
    \begin{aligned}
        \limsup_{\varepsilon \to 0} \lambda_{1,p}(\Omega_\varepsilon) &\leq
         \lim_{\varepsilon \to 0}\frac{\int_0^{R_\varepsilon-\varepsilon} |\dot w(t)|^p \cdot (\delta e^{t-R_\varepsilon})^{n-1} \,\d t}{\int_0^{R_\varepsilon-\varepsilon} |w(t)|^p\cdot (\delta e^{t-R_\varepsilon})^{n-1} \,\d t}
         \\ & =\frac{\int_0^2 |\dot w(t)|^p e^{(n-1)t} \,\d t}{\int_0^2 |w(t)|^p e^{(n-1)t} \,\d t}
         \\ & =\bar\lambda_{2,-1,n}.
    \end{aligned}
    \end{equation}
    Here the last equation is due to Remark \ref{rmk:13} (1). Finally, according to \eqref{eq:thm1.2},
    $$
    \liminf_{\varepsilon \to 0} \lambda_{1,p}(\Omega_\varepsilon) \geq \lim_{\varepsilon \to 0}\bar\lambda_{D_\varepsilon,-1,n}=\bar\lambda_{2,-1,n}
    $$
    since $D_\varepsilon\rightarrow 2$ by \eqref{eq:Diameter}.
    \par Finally, we note that the diameter $2$ used in the above construction can be replaced by any $D > 0.$ Thus the sharpness of the estimate \eqref{eq:thm1.2} is established for all $D > 0.$
    \end{proof}

    \begin{remark}
    \begin{itemize}
      \item [(1)] Recall the classical McKean’s theorem in \cite{McKean1970AnUB}: If $(M,g)$ is an $n-$dimensional complete simply connected Riemannian manifold such that the sectional curvature is bounded above by $K<0,$ then for any bounded domain $\Omega$ in $M,$
          $$
          \lambda_{1,p}(\Omega)\geq\Big(\frac{(n-1)\sqrt{-K}}{p}\Big)^p.
          $$
          The bound is independent of the size or shape of the domain, and depends only on the curvature upper bound $K$ and the dimension $n.$
          \par In contrast, the main theorem of the paper assumes only a Ricci curvature lower bound
$\Ric\geq (n-1)K$ with $K<0,$ and provides a lower bound $\bar\lambda_{D,K,n}.$  This bound decays exponentially as $D\rightarrow\infty,$ reflecting the fact that in the presence of negative Ricci curvature a domain can be made “long and thin” so that its principal frequency becomes arbitrarily small.
\par The sharpness construction in Section 3 exemplifies this difference: the warped product manifolds satisfy
$\Ric\geq-(n-1)$ but have a cusp-like end where sectional curvature is not bounded above by a negative constant (indeed, spherical sectional curvature becomes positive and large). Hence and the eigenvalue can approach zero as the diameter grows.
      \item [(2)] Ona can see Remark 2.1 in \cite{jin2025estimate}. From the same example, we show that the restriction that $\Omega$ is embedded in a non-compact manifold cannot be dropped.
    \end{itemize}

    \end{remark}

\noindent{Xiaoshang Jin}\\
  School of mathematics and statistics, Huazhong University of science and technology, Wuhan, P. R. China. 430074
 \\Email address: jinxs@hust.edu.cn
\bigskip \\
\noindent{Zhiwei L\"u}\\
  School of mathematics and statistics, Huazhong University of science and technology, Wuhan, P. R. China. 430074
 \\Email address: m202470005@hust.edu.cn

\end{document}